# On links of certain semiprime ideals of a noetherian ring

**C.L. Wangneo**

95-C,UNOLane,Talab-Tillo

Jammu,J&K,India,180002.

**Abstract**

In this paper we prove our main theorem, namely, theorem (8), which states that a link $Q \to P$, of prime ideals $Q$ and $P$ of a noetherian ring $R$ that are $\sigma$-semistable with respect to a fixed automorphism $\sigma$ of $R$, induces a link $Q^0 \to P^0$ of the semiprime ideals $Q^0$ and $P^0$ of the ring $R$, where $Q^0$ and $P^0$ are the largest $\sigma$-invariant or $\sigma$-stable ideals contained in the prime ideals $Q$ and $P$. We also prove a converse to this theorem.

**Introduction**

In this paper we study the links of certain semiprime ideals of a noetherian ring R. Following the definition in [1], page 178, or from [3], we recall, there is a link from a prime ideal Q to a prime ideal P of a noetherian ring R ( written as Q→P ), if there is an ideal A of R such that QP≤A<Q∩P and the R- bi-module Q∩P/A is torsion free as a left R/Q - module and as a right R/P - module. Also following [2], we define the σ- invariant or the σ-stable part $I^0$ of an ideal I of R that is σ-semistable with respect to a fixed automorphism σ of R as the largest σ- invariant ideal contained in I. This definition immediately yields the fact that if Q and P are prime ideals of a Noetherian ring R that are σ-semistable with respect to a fixed automorphism σ of R, then $Q^0$ and $P^0$ are σ-invariant semiprime ideals of R. Following the above definition of a link of prime ideals as defined in [1] or [3], we define what we mean by the existence of a link between the semi prime ideals $Q^0$ and $P^0$ of a noetherian ring R, where $Q^0$ and $P^0$ are the σ-invariant parts of the prime ideals Q and P that are assumed to be semistable with respect to an automorphism σ of R. With these definitions in mind we prove our main theorem, namely, theorem (8), which states that if R is a Noetherian ring, then a link Q→P of prime ideals Q and P of R, that are σ-semistable with respect to an automorphism σ of the ring R, induces a link $Q^0 \to P^0$ of the semiprime ideals $Q^0$ and $P^0$ of R. We also prove a converse to this theorem.

**Definitions and Notation**:- To make mention of the source of our reference we state that throughout this paper we adhere to and adapt, as convenience and relevance permits, the notation and definitions of [1],[2], or [3] respectively. Thus, for example, for an ideal I of a noetherian ring R, which is σ-semistable with respect to a fixed automorphism σ of R, we denote by $I^0$, as in [2], lemma(6.9.9), the largest σ- invariant ideal contained in I and call it the σ-invariant part of I. For a right module or a bimodule M over a ring we denote by |M| the right krull dimension of the module M whenever this dimension exists. We must mention here that for a bimodule M over a ring R, we will again use the symbol |M| to denote only the right krull dimension of M unless otherwise stated. For the basic definition and results on krull dimension we refer the reader to [1]. For a few more words about the terminology in this paper we mention that a ring R is noetherian means that R is a left as well as right noetherian ring. If R is a ring and M is a right R module then we denote by Spec.R , the set of prime ideals of R Moreover r-ann.T denotes the right annihilator of a subset T of M and l-ann.T denotes the left annihilator of a left subset T of W in case W is a left R module. For two subsets A and B of a given set, A≤B means B contains A and A<B denotes A≤B but A≠B. For an ideal A of R ,c(A) denotes the set of elements of R that are regular modulo the ideal A. Finally we mention that throughout all our rings are with identity element and all our modules are unitary.

**Main Theorem**

To prove our main theorem we first define the σ-invariant (or the σ-stable) part of an ideal I that is σ-semistable with respect to an automorphism σ of the ring R. Next as mentioned in the introduction, following the definition of a link between two prime ideals of a noetherian ring R as given in [3], for example,we give the definition of a link $Q^o \to P^o$ where $Q^o$ and $P^o$ are the σ - invariant parts of the prime ideals Q and P of a Noetherian ring R which are σ - semi stable with respect to a fixed automorphism σ of R. We then describe when these links exist.

**Definition (1)**: Let R be a noetherian ring and let I be an ideal of R. Then following the definition (6.9.8) of [2], we say I is $\sigma$ - semi-stable or $\sigma$ - semi-invariant ideal of R if there exists an integer $n \geq 1$ such that $\sigma^n(I) = I$. I is said to be $\sigma$ stable or $\sigma$ invariant ideal of R if $\sigma(I) = I$.

With this definition we prove the following result below.

**Proposition (2):** Let R be a Noetherian ring with an automorphism $\sigma$ of R. Then the following hold,

a) If I is a $\sigma$ - semi stable ideal of R with $\sigma^n(I) = I$, $n \geq 1$ and if $I^o = I \cap \sigma(I) \cap \sigma^2(I) \cap \sigma^3(I) \ldots \cap \sigma^{n-1}(I)$, than $I^o$ is the largest $\sigma$-invariant ideal of R contained in I. In case I is a prime ideal of R, then $I^o$ is a $\sigma$ invariant semi prime ideal of R.

b) If I, J are $\sigma$ - semistable ideals of R, then there exists a common integer $k \geq 1$ such that $\sigma^k(I) = I$ and $\sigma^k(J) = J$.

Prof: a) The proof of (a) is obvious (see;e.g;[2],Lemma (6.9.9)(i)).

b) For the proof of (b) note that if m and n integers (m, n $\geq$ 1) such that $\sigma^m(I) = I$ and $\sigma^n(J) = J$, then k = mn implies that $\sigma^k(I) = I$ and $\sigma^k(J) = J$.

**Notation (3):** For an ideal I of R which is $\sigma$-semi stable with respect to an automophism $\sigma$ of R we will denote by $I^o$ the $\sigma$ invariant ideal of proposition (2). We also call Iº the σ invariant part of the ideal I.

**Definition (4):** Let R be a Noetherian ring and let σ be an automorphism of R. Let P, Q $\in$ Spec (R) be prime ideals such that P and Q are σ - semistable prime ideals of R. Let Pº and Qº be the σ - invariant parts of P and Q respectively. Then following the definition of a link between two prime ideals of a noetherian ring R as given, for example, in [3], we say that a link Qº → Pº between the semi prime ideals Qº and Pº exists if there exists a non zero bi-module $Q^0 \cap P^0/A$ with A an ideal of R such that QºPº $\leq$ A < Qº $\cap$ Pº and $Q^0 \cap P^0/A$ is a left $R/Q^0$ torsionfree module and a right $R/P^0$ torsionfree module.

We now prove the result below regarding links of prime ideals of a noetherian ring whose use will become apparent as we proceed to prove our main theorem.

**Proposition (5):** Let R be a Noetherian ring with P, Q $\in$ Spec. R. Let Q → P be a link between the prime ideals Q and P. Then there exists a linking bi-module Q $\cap$ P/B with QP $\leq$ B < Q $\cap$ P such that B is the unique minimal ideal such that a link Q → P exists via the ideal B.

Proof: Consider the set S = {$A_i$/$A_i$ are ideals of R with QP $\leq$ A < Q $\cap$ P and such that a link Q → P exists via the ideal $A_i$}. Let B = $\cap$ $A_i$. Now consider the bi-module, Q $\cap$ P/B.

Clearly QP ≤ B < Q ∩ P. We make the following claim :

Q ∩P/ is a right R/P – torsionfree module.

Proof of the claim: Suppose Q ∩P/B is not right R/P torsionfree module. Then there exists f Є Q ∩ P and f ∉ B such that fg Є B, for some g Є c(P). Now f ∉ B implies that f ∉ Aj for some j. Since B ≤ Aj, so fg Є B implies that fg Є Aj. Also f ∉ Aj and g Є c(P) then means that Q ∩P/ Aj is not a right R/P torsionfree module which contradicts that Aj Є S. hence Q ∩P/B is a right $\frac{R}{P}$ torsionfree module. This proves the claim. Similarly we can show that Q ∩P/B is a left $\frac{R}{Q}$ torsionfree module. Hence Q → P is a link via the ideal B. Moreover it is clear that B is the unique minimal ideal of R such that Q → P is a link via the ideal B because if Q → P is a link via an ideal A of R, then obviously A Є S. which yields that B ≤ A.

**Proposition (6):** Let R be a noetherian ring. Let Q, P Є Spec. R, and let Q → P be a link via the unique minimal ideal A of R. If Q, P are σ- semistable prime ideals of R for an automorphism σ of R, then there exists an integer n ≥1 such that $\sigma^n(p) = P$, $\sigma^n(Q) = Q$ and $\sigma^n(A) = A$.

Proof: First note that because of proposition (2) above there exists a common integer n ≥ 1, such that $\sigma^n(p) = P$ and $\sigma^n(Q) =$

Q. Now the link Q → P via the ideal A implies that for any integer i ≥ 0 there is a link $\sigma^i(Q) \to \sigma^i(P)$ via the ideal $\sigma^i(A)$. In particular $\sigma^n(Q) \to \sigma^n(P)$ via $\sigma^n(A)$ implies that there is a link Q → P via the ideal $\sigma^n(A)$. Hence by the hypothesis on A, A ≤ $\sigma^n(A)$. We now work with the authomorphism $\sigma^{-1}$ and first observe that $\sigma^{-n}(p) = P$ and $\sigma^{-n}(Q) = Q$. The above argument now yields that $A \subseteq \sigma^{-n}(A)$ or $\sigma^n(A) \subseteq A$. Thus we get that $A = \sigma^n(A)$.

**Proposition(7):-** Let R be a Noetherian ring. Let P and Q be prime ideals of R that are σ- semi stable with respect to an automorphism σ of R and let the prime ideal Q be linked to the prime ideal P. Further ( using proposition (5) and proposition (6) above) let Q be linked to P via a unique minimal ideal A and let m be the common integer such that $\sigma^m(Q)=Q, \sigma^m(P)=P$ and $\sigma^m(A)=A$. Let $Q^0, P^0$ and $A^0$ have the usual meaning. Then the following hold true :

(1) The ring R has at most two prime ideals minimal over the ideal A, namely, Q or P. In case P is a minimal prime ideal over A such that |R/A|=|R/P|, then P is also a prime ideal minimal over the ideal $A^0$.

(2) If |R/Q|=|R/P|, then both Q and P are minimal prime ideals over A as well as over the ideal $A^0$. Moreover in this case all

the components of the semiprime ideals $Q^0$ and $P^0$ are minimal prime ideals over $A^0$.

proof:- First observe that it is given that A is the unique minimal ideal of R such that Q → P is a link via the ideal A. Also it is given that Q and P are σ- semistable prime ideals of R, hence, we get that there exists a common integer $m \geq 1$ such that $\sigma^m(Q) = Q$, $\sigma^m(p) = P$ and thus using proposition (6) above we have

$\sigma^m(A) = A$. Let $A^o = A \cap \sigma(A) \cap \ldots \cap \sigma^{m-1}(A)$. We now prove (1).

Proof of (1):- We prove Q or P is a minimal prime ideal over A. To see this we first note that since we are given that Q → P is a link via the ideal A hence we have that QP≤A and this immediately implies that either Q or P is a prime ideal minimal over A. Now suppose that among the prime ideals minimal over A, P is a prime ideal minimal over A, such that |R/P|= |R/A|. Then obviously since |R/A|= |R/A⁰|, hence P is a prime ideal minimal over the ideal $A^0$ as well.

proof of (2):- Now we prove (2) under the assumption that |R/Q|=|R/P|. In this situation then it is clear that either Q=P or Q and P are distinct incomparable prime ideals over the ideal A. Hence, since QP≤A, it is not difficult to see that the set of minimal prime ideals of R/A consists of the

prime ideals Q/A and P/A. Assume that among the prime ideals of R minimal over the ideal A, P is the prime ideal such that $|R/P|=|R/A|$ $(=|R/A^0|)$ then by (1) above P is the prime ideal minimal over the ideal $A^0$ as well. But it is given that $|R/Q|=|R/P|$, Hence again using (1) above we get that Q is also a prime ideal minimal over the ideal $A^0$. The rest is obvious.

We are now ready to prove our main theorem.

**Theorem (8):** Let R be a Noetherian ring. Let P be a prime ideal of R that is σ- semi stable with respect to an automorphism σ of R. Let Q be a σ-semistable prime ideal that is linked to P. Then the link Q → P of the prime ideals Q and P induces a link $Q^o → P^o$ of the semi prime ideals $Q^o$ and $P^o$ of R.

proof: By proposition (5) above we choose the unique minimal ideal A of R such that Q → P is a link via the ideal A. By proposition (2), since Q and P are σ- semi stable prime ideals of R there exists a common integer n ≥ 1 such that

$σ^n(Q) = Q$, $σ^n(p) = P$ and in that case using proposition (6) we have $σ^n(A) = A$. Let $A^o = A \cap σ(A) \cap … \cap σ^{n-1}(A)$. We now prove that $Q^o → P^o$ is a link via the ideal $A^o$.

For the proof we first assume by proposition (7) above that P is a minimal prime ideal over A as well over $A^0$.

We now make the claim:- $A^o \neq Q^o \cap P^o$.

Proof of the claim:-For the proof of this claim we assume that $A^o = Q^o \cap P^o$. Then $A^o = Q^o \cap P^o \leq A < Q \cap P$, implies, that in the semi prime ring $R/A^0$, the prime ideal $Q/A^0$ is linked to the minimal prime ideal $P/A^0$ via the ideal $A/A^0$. But this contradicts lemma (11.17) of [ 1 ]. This proves the claim. Hence we must have $A^o < Q^o \cap P^o$.

We now show that $Q^o \to P^o$ is a link via the ideal $A^o$. But first observe that $Q \to P$ is a link via the ideal A implies also that $\sigma^i(Q) \to \sigma^i(P)$ is a link via the ideal $\sigma^i(A)$ for any integer $i \geq 1$. We now make the claim: $(Q^0 \cap P^0)/A^0$ is a right $R/P^o$ torsionfree module and a left $R/Q^o$ torsionfree module.

Proof of the claim: Suppose first that $(Q^0 \cap P^0)/A^0$ is not a right $R/P^o$ torsionfree module. Then there exists $f \in Q^o \cap P^o$, $f \notin A^o$ and a $g \in c(P^o)$ such that $fg \in A^o$. Now $f \notin A^o$ implies that there exists an integer $m \geq 1$ such that $f \notin \sigma^m(A)$. However, $f \in Q^o \cap P^o$ implies that $f \in \sigma^m(Q) \cap \sigma^m(P)$ and $fg \in A^o$ implies that $fg \in \sigma^m(A)$. Observe that $g \in c(P^o)$ means that $g \in C[\sigma^m(P)]$. All this means that $\sigma^m(Q) \cap \sigma^m(P)/\sigma^m(A)$ is not right $R/\sigma^m(P)$ torsionfree module contradicting our earlier observation that $\sigma^m(Q) \to \sigma^m(P)$ is a right link via the ideal $\sigma^m(A)$. Hence we must have that $(Q^0 \cap P^0)/A^0$ is a right $R/P^o$ torsionfree module. Similarly we can show that $(Q^0 \cap P^0)/A^0$ is a left $R/Q^o$ torsionfree module. Hence $Q^o \to P^o$ is a link via the ideal $A^o$.

There is a converse to the above result which we shall state and prove now. We mention at the outset of this theorem that we shall adapt and mimic, the proof of theorem(11.2)of [1] and hence while doing this we will adapt(without further mention)as much as possible the terminology that is used in the proof of theorem (11.2) of [1].

**Theorem (9):** Let R be a Noetherian ring and let $\sigma$ be an automorphism of R. Let Q, P be $\sigma$ - semi stable prime ideals of R. Then a (right) link $Q^o \to P^o$ of semi prime ideal $Q^o$ and $P^o$ of R implies that for any integer $i \geq o$ there exists an integer $j \geq 0$ such that $\sigma^i(Q) \to \sigma^j(P)$ is a (right) link.

Proof: As stated earlier we will give a sketch of the proof which is on the same lines as the proof of theorem (11.2) of [1].
Since we are given $Q^o \to P^o$ is a (right) link via an ideal A ,thus following [1], theorem (11.2), we assume without loss of generality that A=0. So we may assume $Q^o P^o = 0$ and $Q^o \cap P^o$ is a nonzero torsionfree right $R/P^o$-module and a torsionfree left $R/Q^o$-module. Since l-ann.$(Q^o \cap P^o) = Q^o$, we conclude that l-ann$(Q^o) \leq Q^o$.

Note that l-ann.$(P^o) = Q^o$ because $Q^o P^o = 0$. We first show that $Q^o$ is essential as a right ideal of R. To see this suppose I is a nonzero right ideal of R. If $I Q^o = 0$, then $I <$ l-ann.$(Q^o) \leq Q^o$. Hence $I \cap Q^o \neq 0$. Now clearly if $IQ^o \neq 0$, then obviously $I \cap Q^o \neq 0$. Hence $Q^o$ is essential as a right ideal of R. Next note that $Q^o \cap P^o$ is a torsionfree right $R/P^o$-module and since $Q^o/Q^o \cap P^o$ is isomorphic to a right ideal of $R/P^o$, thus $Q^o/Q^o \cap P^o$ is torsionfree as a right $R/P^o$-module hence $Q^o$ is torsionfree as a right $R/P^o$-module. Thus by proposition (6.18) of [1], $Q^o$ has an essential submodule isomorphic to a finite direct sum of uniform right ideals of $R/P^o$. Since $Q^o_R$ is essential as a right submodule of $R_R$, thus $E(Q^o)_R \approx E(R)_R$, where $E(Q^o)_R$ and $E(R)_R$ are the injective hulls of $Q^o$ and R as right R modules. Now since the components of the semiprime ring $R/P^o$ are isomorphic to each other, we must have that

$E(R)_R \approx E^n$, where E is the injective hull of a uniform right ideal of $R/P^o$ and $n = \text{rank}(R)_R$. Following the argument exactly as in proposition (6.23) of [1], we note that E is independent of the choice of a uniform right ideal of $R/P^o$. Obviously E has an essential submodule which is a torsionfree $R/P^o$-module. Thus ann.$_E(P^o)$ is torsionfree as a $R/P^o$ module. Since l-ann $_R(P^o) = Q^0$, as seen in the first para above, it follows that $Q^o = R \cap \text{ann.}_{E(R)}(P^o)$. Hence, if $W_R = \text{ann.}_{E(R)}(P^o)$, then $R/Q^o = R/R \cap W \approx W + R/W \leq E(R)/W \approx E^n/\text{ann.}_E^n(P^o)$.

Now observe that $E^n/\text{ann.}_E^n(P^o) \approx (E/\text{ann.}_E(P^o))^n$. Hence $R/Q^o$ embeds in $(E/\text{ann.}_E(P^o))^n$. It thus follows that any uniform right ideal of $R/P^o$ embeds in $E/\text{ann.}_E(P^o)$. Now, let K be a submodule of E such that $\text{ann.}_E(P^o) < K$ and $K/\text{ann.}_E(P^o)$ is isomorphic to a uniform right ideal of $R/Q^o$. Choose an element x in K not annihilated by $P^o$. Let $M = xR$ and let $U = \text{ann.}_M(P^o)$. Clearly $M+U/U$ is isomorphic to a uniform right ideal of $R/Q^0$ and by an argument similar to cor.(6.20) of [1] U is isomorphic to a uniform right ideal of $R/P^o$. Since E is uniform, so is M, and it is clear from the definition of U that $0 < U < M$ is an affiliated series for M. As a consequence there exist integers m and n such that U is $R/\sigma^m(P)$-torsionfree and M/U is $R/\sigma^n(Q)$-torsionfree modules. Again by an argument similar to theorem (11.2) of [1], we get that $\sigma^n(Q) \to \sigma^m(P)$. Now applying $\sigma$ repeatedly to the link $\sigma^n(Q) \to \sigma^m(P)$ and observing that Q is $\sigma$-semistable, we get that, for any integer $i \geq 1$ there is a link $\sigma^i(Q) \to \sigma^j(P)$ for some $j \geq 1$. This completes the proof of the theorem.

Remark:- We remark that theorem(11.2) of [1] that we have used in the proof of theorem(9) above actually is proved in [1] to give a characterization of the existence of links between two prime ideals in a noetherian ring and in fact this theorem tells how links really arise in the study of finitely generated modules over a noetherian ring.In this context, thus, our proof of theorem (9) given above is circular. We mention that a direct and a much easier proof of theorem (9),namely one, which does not use theorem(11.2) of [1],may be possible. But we have not verified the same.